\documentclass{amsart}
\usepackage{amssymb}
\usepackage{amsmath}
\usepackage{amssymb}\usepackage{amsmath}\usepackage{a4}
\usepackage{color}

\begin{document}
\date{Version 1r as of 20 February 2008 per PBG}
\newtheorem{theorem}{Theorem}[section]
\newtheorem{lemma}[theorem]{Lemma}
\newtheorem{remark}[theorem]{Remark}
\newtheorem{question}[theorem]{Question}
\def\qedbox{\hbox{$\rlap{$\sqcap$}\sqcup$}}
\makeatletter
  \renewcommand{\theequation}{%
   \thesection.\alph{equation}}
  \@addtoreset{equation}{section}
 \makeatother
\def\Op{\operatorname{Op}}
\title[Heat Content and Isospectrality]
{Heat Content, Heat Trace, and Isospectrality}
\author{P. Gilkey}
\begin{address}{Mathematics Department, University of Oregon, Eugene, OR 97403, USA}\end{address}
\begin{email}{gilkey@darkwing.uoregon.edu}\end{email}
\begin{abstract} We study the heat content function, the heat trace function, and questions of
isospectrality for the Laplacian with Dirichlet boundary conditions on a compact manifold with smooth
boundary in the context of finite coverings and warped products.
\end{abstract}
\keywords{asymptotic expansion, heat content, heat trace, isospectrality, covering projection, warped product.
\newline 2000 {\it Mathematics Subject Classification.} 58J35}
\maketitle
\section{Introduction}
\subsection{The spectral resolution} Let $\mathcal{M}:=(M,g)$ be a compact Riemannian manifold of dimension
$m$ with smooth non-empty boundary
$\partial M$. Let $\operatorname{dvol}_{\mathcal{M}}$ and $\operatorname{dvol}_{\partial\mathcal{M}}$ be the Riemannian measures on $M$ and
on
$\partial M$, respectively. Let $\Delta_{\mathcal{M}}:=\delta d$ be the scalar Laplacian with {\it
Dirichlet boundary conditions}, i.e.
$$\text{Domain}(\Delta_{\mathcal{M}})=\{\phi\in C^\infty(M):\phi|_{\partial M}=0\}\,.$$
There is a complete orthonormal
basis
$\{\phi_n\}$ for
$L^2(\mathcal{M})$ where $\phi_n\in C^\infty(M)$, where
$\phi_n|_{\partial M}=0$, and where
$\Delta_{\mathcal{M}}\phi_n=\lambda_n\phi_n$; these are the {\it Dirichlet eigenfunctions}. The collection
$$\mathcal{S}(\Delta_{\mathcal{M}}):=\{\phi_n,\lambda_n\}$$
is called a {\it spectral resolution} of $\Delta_{\mathcal{M}}$. If one orders the eigenfunctions so
$$0\le\lambda_1\le...\le\lambda_n\ldots,$$
 then one has the Weyl estimate \cite{We15} that $\lambda_n\sim n^{2/m}$ as
$n\rightarrow\infty$. We set
$$\operatorname{Spec}(\Delta_{\mathcal{M}}):=\{\lambda_1,\lambda_2,...\}$$
where each eigenvalue is repeated according to multiplicity.
Two Riemannian manifolds $\mathcal{M}_1$ and $\mathcal{M}_2$ are said to be {\it isospectral} if
$\operatorname{Spec}(\Delta_{\mathcal{M}_1})=\operatorname{Spec}(\Delta_{\mathcal{M}_2})$. We
refer to \cite{GPS05} for further details concerning isospectrality.

\subsection{Operators of Laplace type}
It is convenient to work in slightly greater generality -- this will be important
in Section \ref{sect-3} when we discuss warped products. An operator
$D$ is said to be of {\it Laplace type} if the leading symbol of
$D$ is given by the metric tensor or, equivalently, if we may express in any system of
local coordinates
$x=(x_1,...,x_m)$ the operator $D$ in the form:
$$D=-\{g^{ij}\partial_{x_i}\partial_{x_j}+a^i\partial_{x_i}+b\}$$
where we adopt the {\it Einstein convention} and sum over repeated indices; here $a^i$ and $b$ are smooth
functions and $g^{ij}$ is the inverse of the metric $g_{ij}:=g(\partial_{x_i},\partial_{x_j})$. Let
$$\operatorname{dvol}_{\mathcal{M}}=gdx^1...dx^m$$
where 
$g=\det(g_{ij})^{1/2}$. The scalar Laplacian
$\Delta_{\mathcal{M}}$ is of Laplace type since
\begin{equation}\label{eqn-1.a}
\Delta_{\mathcal{M}}=-g^{-1}\partial_{x_i}gg^{ij}\partial_{x_j}
=-\left(g^{ij}\partial_{x_i}\partial_{x_j}+\{g^{-1}\partial_{x_i}(gg^{ij})\}\partial_{x_j}\right)\,.
\end{equation}

\subsection{The heat equation}
Let $\phi\in C^\infty(M)$ define the
initial temperature of the manifold. Let
$D$ be an operator of Laplace type on $C^\infty(M)$. The subsequent temperature distribution
$u:=e^{-tD}\phi$ for $t>0$ is defined by the following relations, we refer to
\cite{CaJa86} for a further discussion of the heat process:
\begin{equation}\label{eqn-1.b}\begin{array}{ll}
(\partial_t+D)u=0&\text{(evolution equation)},\\
\lim_{t\downarrow0}u(\cdot,t)=\phi\quad\text{in}\quad L^2&\text{(initial condition)},\vphantom{\vrule height 12pt}\\
u(\cdot,t)|_{\partial M}=0&\text{(boundary condition)}\,.\vphantom{\vrule height 12pt}
\end{array}\end{equation}

The special case that $D=\Delta_{\mathcal{M}}$ is of particular interest. Let
$$\sigma_n(\phi):=\int_M\phi(x)\bar\phi_n(x)\operatorname{dvol}_{\mathcal{M}}$$
be the {\it Fourier coefficients}. We may then express
$$
u(x,t)=\sum_ne^{-t\lambda_n}\sigma_n(\phi)\phi_n(x)\,.$$

\subsection{The heat content}
Let $\rho$ be the specific heat and let $\phi$ be the initial temperature of the manifold. The total heat energy content is then defined to be:
\begin{equation}\label{eqn-1.c}
\beta(\phi,\rho,D)(t):=\int_Mu(x;t)\rho(x)\operatorname{dvol}_{\mathcal{M}}\,.
\end{equation}
The heat content is expressible for the Laplacian in terms of the Fourier
coefficients:
$$
\beta(\phi,\rho,\Delta_{\mathcal{M}})(t)=\sum_ne^{-t\lambda_n}\sigma_n(\phi)\sigma_n(\rho)\,.
$$
We shall assume
$\rho$ and $\phi$ are smooth henceforth. We refer to \cite{BeGiSe08} for some results in the non-smooth
setting where
$\phi$ is allowed to blow up near the boundary and to \cite{BeSr90} where the boundary is polygonal. 

The total heat energy $\beta_{\mathcal{M}}(t)$ of $M$ is defined by taking $\phi(x)=\rho(x)=1$;
\begin{equation}\label{eqn-1.d}
\beta_{\mathcal{M}}(t):=\beta(1,1,\Delta_{\mathcal{M}})(t)=\sum_ne^{-t\lambda_n}\sigma_n(1)^2\,.
\end{equation}
The total heat energy content of the manifold is a scalar function which is an isometry invariant of the
manifold.
For example, if $\mathcal{M}=([0,\pi],dx^2)$ is the interval with the standard metric, then 
$$\begin{array}{ll}
\Delta_{\mathcal{M}}=-\partial_x^2,&\displaystyle
\operatorname{Spec}(\Delta_{\mathcal{M}})=\left\{n^2\right\}_{n=1}^\infty,\\
\mathcal{S}(\Delta_{\mathcal{M}})=\left\{\sqrt{\frac2\pi}\sin(nx),n^2\right\}_{n=1}^\infty,&
\sigma_n(1)=\left\{\begin{array}{lll}
2\sqrt{\frac2\pi}&\text{if}&n\equiv1\text{ mod }2\\
0&\text{if}&n\equiv0\text{ mod }2\end{array}\right\},\\
\beta_{\mathcal{M}}(t)=\displaystyle{\textstyle\frac8\pi}\sum_{k=0}^\infty
{\frac{1}{(1+2k)^2}}e^{-(1+2k)^2t}\,.
\end{array}$$

\subsection{The heat trace} Let $D$ be an operator of Laplace type on $C^\infty(M)$. The operator
$e^{-tD}$ is an infinitely smoothing operator. If $f\in C^\infty(M)$ is an auxiliary function which
is used for localization or smoothing, then $fe^{-tD}$ is of trace class and
$\operatorname{Tr}_{L^2}(fe^{-tD})$
is well defined. We shall assume that $f$ is smooth and refer to \cite{BeGiKiSe08} for some results in the
non-smooth setting where
$f$ is allowed to blow up near the boundary. We also refer to \cite{McSi67} for results concerning Riemann
surfaces with corners.

If we take $f=1$ and let $D=\Delta_{\mathcal{M}}$, then
$$\operatorname{Tr}_{L^2}(e^{-t\Delta_{\mathcal{M}}})=\sum_ne^{-t\lambda_n}$$
is a spectral invariant which determines 
$\operatorname{Spec}(\Delta_{\mathcal{M}})$.

\subsection{Local invariants} We can extract locally computable invariants from the heat content and from the
heat trace as follows. Let $D$ be an operator of Laplace type on $C^\infty(M)$ and let $f,\rho,\phi\in
C^\infty(M)$. Work of Greiner
\cite{Gr68} and of Seeley
\cite{Se69, Se69a} can be used to show that there is a complete asymptotic expansion
\begin{equation}\label{eqn-1.e}
\operatorname{Tr}(fe^{-tD})\sim\sum_{n=0}^\infty a_n(f,D)t^{(n-m)/2}
\quad\text{as}\quad t\downarrow0\,.
\end{equation}
Similarly, see the discussion in
\cite{BeGiSe08, BeGeKiKo07}, there is a complete asymptotic expansion:
$$
\beta(\phi,\rho,D)(t)\sim\sum_{n=0}^\infty\beta_n(\phi,\rho,D)t^{n/2}\,.
$$

These invariants are locally computable and have been studied by many authors;
we refer to
\cite{G04} for a more complete discussion of the history of the subject.
To simplify the discussion, we shall only consider the special case
where
$D=\Delta_{\mathcal{M}}$ and where $f=\rho=\phi=1$. We define the following local isometry invariants of the
manifold:
$$a_n(\mathcal{M}):=a_n(1,\Delta_{\mathcal{M}})\quad\text{and}\quad
\beta_n(\mathcal{M}):=\beta_n(1,1,\Delta_{\mathcal{M}})\,.
$$

Let indices $i$, $j$, $k$, $l$ range from $1$ to $m$ and index a local orthonormal
frame $\{e_1,...,e_m\}$ for the tangent bundle of $M$. Let
$R_{ijkl}$ be the components of the Riemann curvature tensor; our sign convention is chosen so that
$R_{1221}=+1$ on the sphere of radius $1$ in $\mathbb{R}^3$. Let $\rho_{ij}:=R_{ikkj}$ be the  {\it Ricci
tensor} and let
$\tau:=\rho_{ii}$ be the {\it scalar curvature}.  Near the boundary we normalize the choice of the local
frame by requiring that $e_m$ is the inward unit geodesic normal. We let indices $a$, $b$, $c$, $d$ range from
$1$ through $m-1$ and index the restricted orthonormal frame $\{e_1,...,e_{m-1}\}$ for the tangent bundle of
the boundary. Let $L_{ab}:=g(\nabla_{e_a}e_b,e_m)$ be the components of the second fundamental form. We can
use the Levi-Civita connection on $\mathcal{M}$ to multiply covariantly differentiate a tensor defined in the
interior; we let `;' denote the components of such a tensor. Similarly, we can use the Levi-Civita connection
of $\partial\mathcal{M}:=(\partial M,g|_{\partial M})$ to multiply covariantly differentiate a
tensor defined on the boundary; we let `:' denote the components of such a tensor. The difference between `;'
and `:' is measured by the second fundamental form.

\begin{theorem}\ \label{thm-1.1}
\begin{enumerate}
\item $a_0(\mathcal{M})=(4\pi)^{-m/2}\operatorname{Volume}(M)$.
\smallbreak\item 
$a_1(\mathcal{M})=-\frac14(4\pi)^{-(m-1)/2}\operatorname{Volume}(\partial M)$. 
\smallbreak\item $a_2(\mathcal{M})=\frac16(4\pi)^{-m/2}\left\{\int_M\tau 
\operatorname{dvol}_{\mathcal{M}}
+\int_{\partial M}2L_{aa}\operatorname{dvol}_{\partial\mathcal{M}}\right\}$.
\smallbreak\item $a_3(\mathcal{M})=-\frac1{384}(4\pi)^{-(m-1)/2}
\int_{\partial M}\{16\tau+8R_{amam}
+7L_{aa}L_{bb}$\smallbreak$-10L_{ab}L_{ab}\}\operatorname{dvol}_{\partial\mathcal{M}}$.
\smallbreak\item
$a_4(\mathcal{M})=\frac1{360}(4\pi)^{-m/2}\int_M\{12\tau_{;kk}+5\tau^{2}-2|\rho^{2}|+2|R^{2}|\}
\operatorname{dvol}_{\mathcal{M}}$
\smallbreak
$+\frac1{360}(4\pi)^{-m/2}\int_{\partial M}\{-18\tau_{;m}+20\tau L_{aa}
+4R_{am
am}L_{bb}$\smallbreak$ -12R_{am bm}L_{ab}+4R_{ab
cb}L_{ac}+24L_{aa:bb}+\frac{40}{21}L_{aa}L_{bb}L_{cc}$\smallbreak$
-\frac{88}7L_{ab}L_{ab}L_{cc}+\frac{320}{21}L_{ab}L_{bc}L_{ac})\}\operatorname{dvol}_{\partial\mathcal{M}}$.
\end{enumerate}
\end{theorem}
\begin{theorem}\ \label{thm-1.2}
\begin{enumerate}
\smallbreak\item $\beta_0(\mathcal{M})=\operatorname{Volume}(M)$.
\smallbreak\item $\beta_1(\mathcal{M})=
 -\frac2{\sqrt\pi}\operatorname{Volume}(\partial M)$.
\smallskip\item $\beta_2(\mathcal{M})={\int}_{\partial M}\frac12
 L_{aa}\operatorname{dvol}_{\partial\mathcal{M}}$.
\smallbreak\item $\beta_3(\mathcal{M})=
 -\frac2{\sqrt\pi}{\int}_{\partial M}\{\frac1{12} L_{aa}L_{bb}
 -\frac16 L_{ab}L_{ab}-\frac16 R_{amma}\}\operatorname{dvol}_{\partial\mathcal{M}}$.
\smallbreak\item $\beta_4(\mathcal{M})={\int}_{{\partial M}}\{- \frac1{16} L_{
ab}L_{ab}L_{cc}+\frac18 L_{ab}L_{ac}L_{bc}-\frac1{16} R_{am bm}L_{
ab}$\smallbreak$
+\frac1{16} R_{ab cb}L_{ac}+\frac1{32}\tau_{;m}\}\operatorname{dvol}_{\partial\mathcal{M}}$.
\end{enumerate}\end{theorem}

Although formulas for $a_5(\mathcal{M})$ and $\beta_5(\mathcal{M})$ are known, we have
omitted them in the interests of brevity.  Formulas generalizing those in Theorems
\ref{thm-1.1} and
\ref{thm-1.2} are available in the more general setting to study the invariants
$a_n(f,D)$ and $\beta_n(\phi,\rho,D)$ for an arbitrary vector valued operator $D$ of Laplace type; again, we
shall omit details in the interests of brevity and instead refer to the discussion in
\cite{G04}. Although we have chosen to work with Dirichlet boundary conditions, similar formulas exist for
Neumann, transfer, transmittal, and spectral boundary conditions. The  history of this subject is a vast one
and beyond the scope of this brief article to give in any depth. We refer to \cite{Gru96} for a more detailed
discussion of elliptic boundary conditions.

\subsection{Relating the heat trace and heat content}
McDonald and Meyers \cite{MM03a} have constructed additional
invariants involving exit time moments which determine both the heat trace and the heat content; we also refer
to related work
\cite{MM03} by these authors in the context of graphs. 

It is difficult in general, however, to relate the
heat trace and the heat content directly. In particular, there is no obvious relation between the formulas
given in Theorems \ref{thm-1.1} and \ref{thm-1.2} when $n\ge3$. It is clear that
$\operatorname{Tr}\{e^{-t\Delta_{\mathcal{M}}}\}$ is determined by
$\operatorname{Spec}(\Delta_{\mathcal{M}})$ and it is clear that
$\beta_{\mathcal{M}}(t)$ is determined by the full spectral resolution $\mathcal{S}(\Delta_{\mathcal{M}})$. It is not known, 
however, if the full heat content function $\beta_{\mathcal{M}}(t)$ or in particular the
heat content asymptotic coefficients
$\beta_k(\mathcal{M})$ might be determined by $\operatorname{Spec}(\Delta_{\mathcal{M}})$ alone. More specifically, one
does not know if there are Dirichlet isospectral manifolds with different heat content functions.
In the remainder of this brief note, we shall present some results which relate to this question. In Section
\ref{sect-2} we discuss finite coverings and in Section \ref{sect-3} we discuss warped products.

\section{Finite coverings}\label{sect-2}
\subsection{Notational conventions}
We
suppose that
$\pi:M_1\rightarrow M_2$ is a finite
$k$-sheeted covering of compact manifolds with boundary. We assume that $M_2$ is equipped with a Riemannian
metric $g_2$ and choose the induced metric $g_1:=\pi^*g_2$ on $M_1$. Thus $\pi$ is a local isometry and
$\operatorname{Volume}(M_1)=k\operatorname{Volume}(M_2)$.  Since 
\begin{equation}\label{eqn-2.a}
|\pi^*\phi|^2_{L^2(\mathcal{M}_1)}=k|\phi|^2_{L^2(\mathcal{M}_2)},
\end{equation}
pullback $\pi^*$ is an injective closed map from $L^2(\mathcal{M}_2)$ to $L^2(\mathcal{M}_1)$. 

\subsection{Heat trace and heat content asymptotics} The invariants $a_n(\mathcal{M})$ and
$\beta_n(\mathcal{M})$ are locally computable. Since integration is multiplicative under finite coverings,
the following result is immediate:

\begin{theorem}\label{thm-2.1} Let $\mathcal{M}_1\rightarrow\mathcal{M}_2$ be a finite $k$-sheeted Riemannian
cover. Then
$a_n(\mathcal{M}_1)=ka_n(\mathcal{M}_2)$ and $\beta_n(\mathcal{M}_1)=k\beta_n(\mathcal{M}_2)$ for all
$n$.
\end{theorem}

\subsection{Heat trace}\label{sect-2.3}
We begin by presenting an example to show that there are
examples where $\operatorname{Tr}_{L^2}(e^{-t\Delta_{M_1}})\ne k\operatorname{Tr}_{L^2}(e^{-t\Delta_{M_2}})$
despite the fact that the heat content function is multiplicative under finite coverings. Let 
$$\mathcal{M}_1:=([0,4\pi],d\theta^2]/0\sim4\pi
\quad\text{and}\quad
\mathcal{M}_2:=([0,2\pi],d\theta^2)/0\sim2\pi;$$
$\mathcal{M}_1$ may be identified with the circle of radius $2$ in $\mathbb{R}^2$ and $\mathcal{M}_2$
may be indentified with the circle of $1$ in $\mathbb{R}^2$. The natural projection from
$\mathcal{M}_1\rightarrow\mathcal{M}_2$ can be regarded as the double cover of the circle by the circle
induced by the map $z\rightarrow\frac14z^2$. Then
$$\begin{array}{ll}
\textstyle\mathcal{S}(\Delta_{\mathcal{M}_1})=
   \left\{\frac1{\sqrt{4\pi}}e^{\sqrt{-1}k\theta/2},\frac14k^2\right\}_{k=-\infty}^\infty,&
\displaystyle  \operatorname{Tr}_{L^2}\left\{e^{-t\Delta_{\mathcal{M}_1}}\right\}=1+2\sum_{k=1}^\infty
e^{-\frac14k^2t},\\
\textstyle\mathcal{S}(\Delta_{\mathcal{M}_2})=
   \left\{\frac1{\sqrt{2\pi}}e^{\sqrt{-1}k\theta},k^2\right\}_{k=-\infty}^\infty,&
\displaystyle
\operatorname{Tr}_{L^2}\left\{e^{-t\Delta_{\mathcal{M}_2}}\right\}=1+2\sum_{k=1}^\infty e^{-k^2t}\,.
\end{array}$$
Consequently 
$$\operatorname{Tr}_{L^2}\left\{e^{-t\Delta_{\mathcal{M}_1}}\right\}\ne
2\operatorname{Tr}_{L^2}\left\{e^{-t\Delta_{\mathcal{M}_1}}\right\}\,.
$$

Although this example is in the category of closed
manifolds, we can construct other examples as follows. Let $\mathcal{N}=([0,\pi],d\theta^2)$ be a manifold
with boundary. Let
$\tilde{\mathcal{M}}_i:=\mathcal{N}\times\mathcal{M}_i$ and let $\pi$ act only on the second factor. Since
$$e^{-t\Delta_{\mathcal{N}\times\mathcal{M}}}=e^{-t\Delta_{\mathcal{N}}}e^{-t\Delta_{\mathcal{M}}}\,,$$
one has:
\begin{eqnarray*}
&&\operatorname{Tr}_{L^2}\left\{e^{-t\Delta_{\tilde{\mathcal{M}_1}}}\right\}
=\sum_{\ell=1}^\infty e^{-tk^2}\cdot\left\{1+2\sum_{k=1}^\infty e^{-\frac14k^2t}\right\}\\
&\ne&
2\sum_{\ell=1}^\infty e^{-tk^2}\cdot\left\{1+2\sum_{k=0}^\infty e^{-k^2t}\right\}
=2\operatorname{Tr}_{L^2}\left\{e^{-t\Delta_{\tilde{\mathcal{M}_2}}}\right\}\,.
\end{eqnarray*}

\subsection{Heat content asymptotics}
It is perhaps somewhat surprising that in contrast to the situation with the heat trace asymptotics discussed
in Section \ref{sect-2.3} that one has:
\begin{theorem}\label{thm-2.2}
 Let $\mathcal{M}_1\rightarrow\mathcal{M}_2$ be a finite $k$-sheeted Riemannian
cover. Then
$\beta(\mathcal{M}_1)(t)=k\beta(\mathcal{M}_2)(t)$.
\end{theorem}

\begin{proof} 
Let
$\{\lambda_n,\phi_n\}$ be a spectral resolution of
$\Delta_{\mathcal{M}_2}$. Let $c_n=\sigma_n(1)$ be the associated Fourier coefficients. We use
Equation (\ref{eqn-2.a}) to see that
$$1=\sum_nc_n\phi_n\quad\text{in}\quad
L^2(\mathcal{M}_2)\quad\text{implies}\quad1=\sum_nc_n\pi^*\phi_n\quad\text{in}\quad
L^2(\mathcal{M}_1)\,.$$
Since
$\Delta_{\mathcal{M}_2}\pi^*\phi_n=\pi^*\Delta_{\mathcal{M}_1}\phi_n=\lambda_n\pi^*\phi_n$ and since
$\pi^*\phi_n$ satisfy Dirichlet boundary conditions, we
have
$$\{e^{-t\Delta_{\mathcal{M}_1}}\}1=\sum_ne^{-t\lambda_n}c_n\pi^*\phi_n=\pi^*\{e^{-t\Delta_{\mathcal{M}_2}}\}1\,.$$
Consequently
\medbreak\qquad
$\beta_{\mathcal{M}_1}(t)=\displaystyle\langle
e^{-t\Delta_{\mathcal{M}_1}}\pi^*1,\pi^*1\rangle_{L^2(\mathcal{M}_1)} =\langle
\pi^*e^{-t\Delta_{\mathcal{M}_2}}1,\pi^*1\rangle_{L^2(\mathcal{M}_1)}$\par\medbreak\qquad\qquad
\phantom{...}
$=k\langle e^{-t\Delta_{\mathcal{M}_2}}1,1\rangle_{L^2(\mathcal{M}_2)}=k\beta(\mathcal{M}_2)(t)$.\hfill
\end{proof}

\subsection{Summary} Theorems \ref{thm-2.1} and \ref{thm-2.2} show that a Sunada construction involving finite
coverings will not produce isospectral manifolds with different heat content functions as only the order of
the cover is detected. If $\mathcal{M}$ is a Riemannian manifold which has constant sectional curvature $+1$,
then
$\mathcal{M}$ is said to be a {\it sperical space form}. If $\mathcal{M}$ is closed and if the fundamental
group $\pi_1(M)$ is cyclic, then
$\mathcal{M}$ is said to be a {\it lens space}. Ikeda \cite{Ik80,Ik83} and other authors have studied
questions of isospectrality for spherical space forms; we refer to \cite{G95, G04} for further details as the
literature is an extensive one. These examples can easily be modified to the category of manifolds with
boundary by punching out a small disk from
$M_2$ and then lifting to get a spherical spaceform with boundary. Since there are spherical space forms with
the same fundamental group which are not isospectral, neither the heat trace asymptotics nor the full heat
content function determine either the spectrum of the manifold or the isometry type of the manifold.

\section{Warped product metrics}\label{sect-3}
\subsection{Notational conventions}
Let
$\mathcal{N}=(N,g_N)$ be a smooth Riemannian manifold  of dimension $n$ with smooth boundary $\partial N$, let
$\mathcal{M}=(M,g_M)$ be a closed Riemannian manifold of dimension $m$, and let $f\in C^\infty(N)$. We
consider the warped product
$$\mathcal{N}\times_f\mathcal{M}:=(N\times M,g_N+e^{\frac2mf}g_M)\,.$$
The normalizing constant $\frac2m$ is chosen so that one has the following relationship between the volume
elements:
\begin{equation}\label{eqn-3.a}
\operatorname{dvol}_{\mathcal{N}\times_f\mathcal{M}}=e^f\operatorname{dvol}_{\mathcal{N}}
\cdot\operatorname{dvol}_{\mathcal{M}}
\end{equation}

We define an auxiliary operator of Laplace type on $C^\infty(N)$ by setting:
$$D_{\mathcal{N},f}:=e^{-f}\Delta_{\mathcal{N}}e^{f}\,.$$
Note that this operator is no longer self-adjoint if $f$ is non-constant; this operator does, however, have the
same spectrum as $\Delta_{\mathcal{N}}$ since it is conjugate to this operator. We may then use Equation
(\ref{eqn-1.a}) to see that
\begin{equation}\label{eqn-3.b}
\Delta_{\mathcal{N}\times_f\mathcal{M}}=D_{\mathcal{N},f}+e^{-\frac2mf}\Delta_{\mathcal{M}}\,.
\end{equation}

\subsection{The heat content}
Let $\beta(\phi,\rho,D)(t)$ be the generalized heat content
function defined in Equation (\ref{eqn-1.c}).

\begin{theorem}\label{thm-3.1}\ \begin{enumerate}
\item
$\beta(\mathcal{N}\times_f\mathcal{M})(t)=\operatorname{Volume}(\mathcal{M})\cdot\beta(1,e^{f},D_{N,f})(t)$.
\item If $\operatorname{Volume}(\mathcal{M}_1)=\operatorname{Volume}(\mathcal{M}_1)$, then
$\beta(\mathcal{N}\times_f\mathcal{M}_1)(t)=\beta(\mathcal{N}\times_f\mathcal{M}_2)(t)$.
\end{enumerate}
\end{theorem}

\begin{proof} Let
$$u:=e^{-tD_{\mathcal{N},f}}\cdot1$$ be the solution of Equation (\ref{eqn-1.b}) on
$N$ with initial condition
$\phi(\cdot)=1$ which is defined by the operator $D_{\mathcal{N},f}$. Extend
$u$ to $N\times M$ to be independent of the second variable. We apply Equation (\ref{eqn-3.b}). Since
$\Delta_{\mathcal{M}}u=0$, $u$ satisfies Equation (\ref{eqn-1.b}) on $N\times M$ with
initial condition $\phi(\cdot)=1$ using the operator $\Delta_{\mathcal{N}\times_f\mathcal{M}}$. Thus we also
have that
$$u=e^{-t\Delta_{\mathcal{N}\times_f\mathcal{M}}}\cdot1\,.$$
One may now use Equation (\ref{eqn-3.a}) to compute:
\begin{eqnarray*}
\beta(\mathcal{N}\times_f\mathcal{M})(t)
&=&\int_{N\times M}u(x_N;t)e^f\operatorname{dvol}_{\mathcal{N}}\operatorname{dvol}_{\mathcal{M}}\\
&=&\operatorname{Volume}(M)\int_Nu(x_N;t)e^f\operatorname{dvol}_{\mathcal{N}}\\
&=&\operatorname{Volume}(M)\cdot\beta(1,e^f,D_{\mathcal{N},f})(t)\,.
\end{eqnarray*}
This establishes Assertion (1); Assertion (2) follows from Assertion (1).
\end{proof}

Theorem \ref{thm-3.1} shows that the heat content does not even determine the dimension of the underlying
manifold as only the volume of the manifold $M$ appears in this formula. On the other hand, Equation
(\ref{eqn-1.e}) shows that the dimension of the underlying manifold is determined by the heat trace.
Consequently, we once again see that the heat content function does not determine the underlying spectrum.

\subsection{Isospectrality}
We conclude our discussion by showing that isospectrality is preserved by the warped product construction.

\begin{theorem}\label{thm-3.2}\ 
\begin{enumerate}
\item Let $\displaystyle\operatorname{Spec}(\Delta_{\mathcal{M}})=\left\{\mu_i\right\}_{i=1}^\infty$. Then
 $$\operatorname{Spec}\left(\Delta_{\mathcal{N}\times_f\mathcal{M}}\right)=\bigcup_{i=1}^\infty
\operatorname{Spec}\left(\Delta_{\mathcal{N}}+\mu_ie^{-\frac2mf}\right)\,.$$
\item  If
$\operatorname{Spec}\{\mathcal{M}_1\}=\operatorname{Spec}\{\mathcal{M}_2\}$, then
$$
\operatorname{Spec}\{\Delta_{\mathcal{N}\times_f\mathcal{M}_1}\}=
\operatorname{Spec}\{\Delta_{\mathcal{N}\times_f\mathcal{M}_2}\}\,.
$$
\end{enumerate}
\end{theorem}

\begin{proof} Let $\mathcal{M}$ be a Riemannian
manifold. Let $\{\Phi_i,\mu_i\}$ be a spectral resolution of $\Delta_{\mathcal{M}}$. We decompose
\begin{equation}\label{eqn-3-x3}
L^2(\mathcal{N}\times_f\mathcal{M})=\oplus_i L^2(\mathcal{N})\cdot\Phi_i\,.
\end{equation}
Let $\mu_ie^{-\frac2mf}$ act by scalar multiplication. We use Equation (\ref{eqn-3.b}) to see that the
decomposition of Equation (\ref{eqn-3-x3}) induces a corresponding decomposition
$$\displaystyle\Delta_{\mathcal{N}\times_f\mathcal{M}}=\oplus_i
\left\{e^f\Delta_{\mathcal{N}}e^{-f}+\mu_ie^{-\frac2mf}\right\}\,.$$
Assertion (1) now follows since
$$\operatorname{Spec}\left\{e^f\Delta_{\mathcal{N}}e^{-f}+\mu_ie^{-\frac2mf}\right\}
=\operatorname{Spec}\left\{\Delta_{\mathcal{N}}+\mu_ie^{-\frac2mf}\right\}\,.$$
Assertion (2) follows from Assertion (1).
\end{proof}

We may take $N=[0,\pi]$ and assume that $f(0)=f(\pi)=0$. We then have that $\partial(N\times M)$ is isometric
to the disjoint union of two copies of $M$. Since there are many pairs of isospectral closed manifolds which
are not isometric, Theorem
\ref{thm-3.2} provides examples of isospectral manifolds with boundary given by warped products which are not
isometric.

\subsection{Conclusion} Theorems \ref{thm-1.1} and \ref{thm-1.2} show that the volume of the
interior, the volume of the boundary, and the dimension of $M$ are determined by the heat trace. Thus Theorem
\ref{thm-3.1} shows that a warped product construction involving isospectral manifolds with a suitably chosen
manifold with boundary will not produce isospectral manifolds with different heat content functions. Theorem
\ref{thm-3.1} does show, however, that there exist manifolds with the same heat content function which are not
isospectral. If we take $f=1$ and apply the argument of Theorem \ref{thm-3.1}, we see that the heat content
function does not determine the dimension of the manifold.

There exist spherical space forms $\mathcal{M}_1$ and $\mathcal{M}_2$ which are isospectral but not
diffeomorphic. If we take
$\mathcal{N}=([a,b],dx^2)$ with $0<a<b$ and if we take as a warping function
$f(x)=x^2$, then the resulting warped products $\mathcal{P}_i:=\mathcal{N}\times_f\mathcal{M}_i$ are flat
isospectral manifolds whose boundaries are not diffeomorphic.

\section*{Acknowledgments}Research of P. Gilkey was partially
supported by the Max Planck Institute for Mathematics in the Sciences (Germany) and by Project MTM2006-01432 (Spain).

\bibliographystyle{amsalpha}

\end{document}